\documentclass{article}
\usepackage[utf8]{inputenc}
\usepackage{multirow}

\usepackage{graphicx} 
\usepackage{amsmath}

\usepackage{url}
\usepackage[numbers]{natbib} 

\usepackage{appendix}

\title{Meeting Energy Needs by Balancing Cost and Sustainability through Linear Programming}
\author{Daria Soboleva, Alice Hoover, Melinda Koelling}
\date{}
\begin{document}

\maketitle

\begin{abstract}
    Greenhouse gases (GHG) trap heat and make the planet warmer, exacerbating global climate change. Energy production is the second largest contributor to climate change \cite{abstract_info}. In particular, the production of electricity and use of gas contribute to climate change.  Additionally, gas is not renewable and the source of electricity may not be renewable either.  How and whether communities transition to renewable and/or cleaner energy sources is dependent on a number of factors, including energy needs, cost, space, emissions, jobs, and materials.  In this paper, we explore minimizing the cost of building, operating, and maintaining energy sources.  We consider different combinations of cleaner and/or more renewable energy sources to meet energy needed for a given city while keeping total emissions, land use, and energy infrastructure costs low.  For specificity, we use the city of Chicago as a test case.  If we use a combination of wind and solar energy to meet the total energy needs of Chicago, we find that it is most cost effective to use only wind.  However when variable demand and production are included, it is most cost effective to use a combination of wind and solar.  If nuclear and geothermal energy are included to decrease overproduction, it is most cost effective to use a combination of wind and geothermal energy.  

\end{abstract}

\section{Introduction} \label{section_1_introduction}

Renewable energies have been around for quite some time; the very first commercial sale of wind turbines took place in 1927. However, most energy that is produced in the world comes from non-renewable sources. Specifically,  global energy production is fuelled 31.2$\%$ by oil, 27.2$\%$ by coal, and 24.7$\%$ by natural gas \cite{EnergyProdWorld}. 

Energy and its availability has always been taken for granted \cite{Superorganism}. Access to energy sources sets the physical limits to the society: all life, commerce, work, and technological innovation are possible and limited by the new energy available \cite{Energy_is_essential}, \cite{Superorganism}. Moreover, energy consumption per capita has been increasing steadily over recent decades \cite{Consumption_is_growing}, \cite{Consumption_is_growing_2}, especially in high and upper-middle income countries. On the other side, easily available or existing non-renewable energy sources are depleting \cite{Depletion_1}, \cite{Depletion_2}. The International Energy Agency suggests that with no new drilling, the world's oil production would be reduced to 50$\%$ by 2025 and to 15$\%$ by 2040 \cite{oil_no_new_drillings}. Although new sources of energy will be found, the raw materials will be more expensive to extract \cite{Depletion_2}. This will inevitably affect economic growth, given that we have already found and used most of the cheapest sources. This warrants consideration for renewable energies as partial substitution of the non-renewable sources. Not only are renewable energy sources sustainable over time, but they significantly reduce harmful emissions, which is an important factor given that current emissions from energy production grow each year and have negative effects on the environment \cite{EmissionsOtherSources}.

In this project, our goal is to explore how major cities can meet their electricity and gas usage needs exclusively with renewable or clean sources, while minimizing the operation and maintenance cost and meeting several constraints. Although renewable energy cannot eliminate the need for non-renewable sources \cite{energy_cannot_be_fully_substituted}, we attempt to substitute the electricity and natural gas needs (referred to as ``energy needs" from now on). This excludes other direct or indirect consumption of energy by the population, such as transportation, purchased essentials (such as clothing, food, home goods), and recreation.  As an example, we will be exploring the energy needs of the city of Chicago, Illinois. 

In this paper, we first consider a model where the only goal is to minimize costs, in USD per MWh, mainly to show why non-renewable energy is the primary choice of energy source. We then add constraints to the linear model, including emissions, budget, space, and energy needs. We start with considering exclusively solar and wind systems, and later add nuclear and geothermal to solve issues related to overproduction.

\section{Methods: Linear Programming} \label{section_2_methods}
We will be using linear programming to explore solutions in this paper.
Linear programming is a set of algorithmic procedures utilized for solving optimization problems, which involve the maximization or minimization of a linear objective function subject to linear constraints. Maximization problems often deal with profit, while minimization problems often deal with cost, such as the problem we consider in this paper.
More specifically, we use the simplex algorithm for solving linear programs. The simplex algorithm was created in the 1940's by George B. Dantzig. It produces an optimal solution to a linear program (provided such a solution exists) in a reasonable number of steps and is easily implemented on a computer \cite{Strayer_book}.

To solve the models constructed in this paper, we use Python 2021.2 (community edition), and build our code on the existing modeller ``PuLP 2.7.0", written in Python to solve linear programming problems.

\section{First Model} \label{section_3_very_first_model}

Although we consider emissions an important aspect of this project, the financial side is a practical limitation.
Renewable energy infrastructure can be expensive to build, but it also has operation and maintenance costs.
Our objective function minimizes the construction, operation, and maintenance costs. 
For the purpose of this paper, we will call this total cost.

\subsection{Minimizing Total Cost}
This model includes energy from wind, solar, nuclear, geothermal, natural gas, and hydroelectric sources.  
Variables measure how many megawatt hours (MWh) of energy from each source is used in a year:

\begin{tabular}{cl}
    $x_1=$  &  MWh from wind in a year\\
      $x_2=$  &  MWh from solar in a year\\
       $x_3=$  &  MWh from nuclear in a year\\
        $x_4=$  &  MWh from geothermal in a year\\
         $x_5=$  &  MWh from gas in a year\\
          $x_6=$  &  MWh from hydroelectric in a year.
\end{tabular}

\noindent The levelized cost of energy (LCOE) per MWh accounts for construction, operation, and maintenance of each energy source over the course of its lifetime \cite{Cost_OM_2}. The cost per MWh in 2021 for each source is as follows in Table \ref{Tbl:Costs_energy} \cite{Cost_OM_1}.

{
\begin{table}[ht]

\centering
\begin{tabular}{|p{4cm}|p{5cm}|}
\hline
cost per MWh & source of energy \\ \hline 
37.80  &  wind (onshore)\\ 
58.62  &  solar (hybrid) \\ 
96.2  &  nuclear\\ 
39.61  &  geothermal\\ 
37.50  &   combined cycle \\ 
63.9  & hydroelectric \\
\hline
\end{tabular}
\caption{Emissions from conventional energy sources  \protect\cite{Cost_OM_1}}
\label{Tbl:Costs_energy}
\end{table}
}

\noindent  With the variables and levelized cost above, the total cost is

\begin{equation}
  37.80 x_1 + 58.62 x_2+ 96.2 x_3 + 39.61 x_4 + 37.50 x_5 + 63.9 x_6.
  \label{initial_objective_equation_1}
\end{equation}

\noindent If the total needs of Chicago is E MWh in a year, 
then the variables are constrained by 

\begin{equation}
    x_1+x_2+x_3+x_4+x_5+x_6 \geq E
\end{equation}
with $x_i\geq 0$ for all $i$.

With these constraints, the minimum of total cost occurs when we use exclusively non-renewable and non-clean (combined cycle natural gas) energy. This result makes sense: we currently rely on so much non-renewable, non-clean energy because it is the most economical option at present.

\subsection{Possible additional constraints} 
Since we aim to study reducing emissions, using only a non-clean energy source is not acceptable. In the rest of this study, a constraint about emissions will be included. By including clean energy sources, other concerns arise.  As a result, other models in this paper include some of the following.   
\begin{enumerate}
\item 
Meet Chicago's energy needs. This means that the total amount of MWh produced by renewable energy sources is greater than or equal to the energy needs of Chicago.  
\item 
Keep carbon emissions low. Emissions have a negative impact on the environment, and are a significant motivating factor for the shift toward clean and sustainable energy production.The Table \ref{Table:emissions} shows the average emissions of CO$_2$ per MWh of major non-renewable, non-clean sources, as well as renewable and/or clean sources: 

{
\begin{table}[ht]

\centering
\begin{tabular}{  |p{4cm}||p{6cm}| }
 \hline
 Source of energy & Emissions per MWh (g CO$_2$ per MWh)) \\
 \hline
 Coal    & 820,000   \\
 Natural gas &   490,000  \\
 Oil (petroleum) &   1,106,765  \\
 Biomass &  230,000 \\
 \hline
 Wind   &   4970 \\
 Solar  & 45,000 \\
 Nuclear   & 42,200 \\
 Geothermal   & 38,000 \\
 \hline

\end{tabular}
\caption{Emissions from conventional energy sources \protect\cite{EmissionsOtherSources}, \protect\cite{EmissionsOtherSources_2}.}
\label{Table:emissions}
\end{table}
} 

\item 
Meet a certain budget. This is the amount of money that needs to be invested in the very beginning to build all necessary infrastructure. It does not include operation and maintenance costs.
\item 
Keep infrastructure in available/reasonable space.
\item 
Provide for continuously available energy. Since energy needs change throughout the day, we need to adjust our model accordingly so there is no underproduction.
\end{enumerate}

\section{Using solar and wind systems to meet energy needs} \label{section_4_model_using_solar_and_wind_no_space_constraint}

The model consists of a simplified version of objective function (\ref{initial_objective_equation_1}), which only includes wind (x$_1$, in MWh produced a year) and solar (x$_2$, in MWh produced a year) energy. There are five constraints: one constraint for energy needs of Chicago, one for emissions, one for budget, and two for space. The next subsection \ref{section_4_constraints_development} will cover each constraint in detail. See section \ref{section_4_linear_program_model} for the linear model.

\subsection{Determining Constraints} \label{section_4_constraints_development}

\subsubsection{Chicago energy needs} \label{section_4_constraints_chicago_needs}

This constraint ensures that the annual energy production is greater than or equal to the annual energy needs of Chicago. In order to set the lower bound for the energy needs of Chicago, we estimate how much energy (MWh) is used annually. 

From the State Energy Data Portal (SEDS) database, Illinois consumed 1,091,285,298.314 MWh in 2021 \cite{Illinois_demand_2010}. 
Since this database does not include Chicago's consumption, we make an estimate based on Illinois' energy use. 

Chicago Data Portal provides a database on energy consumption in Chicago from 2010 \cite{Chicago_demand_2010} but nothing more recent. 
Since SEDS provides the the consumption of Illinois in 2010, we calculate the percentage that went to Chicago that year.
The state of Illinois consumed 1,168,009,546 MWh in 2010 \cite{Illinois_demand_2010}.
That same year, buildings in Chicago consumed 15,142,030 MWh of electricity and 37,998,300 MWh of gas \cite{Chicago_demand_2010}. 
Since the electricity data comprises 68 percent of overall electrical usage in the city and gas data comprises 81 percent of all gas consumption in Chicago for 2010, these values are adjusted and summed to find that 
%
%
%
Chicago consumed about 68,771,765 MWh in 2010.
Dividing Chicago's energy use by Illinois', we estimate that approximately 5.88$\%$ of Illinois' total usage of energy in 2010 went to Chicago. This can be found in Table \ref{tab:consumtpion_gas_electricity}.

{
\begin{table}[ht]
\centering
\begin{tabular}{ |p{3cm}||p{2.5cm}|p{2.5cm}|p{2.5cm}|  }
 \hline

 & Consumption reported (MWh) & $\%$ of total consumption & Estimated total (MWh) consumption\\
 \hline
 Electricity   & 15,142,030  & 68$\%$ & 22,267,691\\

 Natural gas &   37,998,300   & 81$\%$ &46,504,074\\
 \hline

\end{tabular}
\caption{Estimated consumption of electricity and gas in 2010, as computed above.}
\label{tab:consumtpion_gas_electricity}
\end{table}
} 




Next, the percentage is adjusted due to population growth.
The population of Chicago has increased by 4$\%$ since 2010 \cite{Chicago_population}, so we raise the percentage from 5.88$\%$ to 7$\%$. 

Applying this to Illinois' energy needs, we find that Chicago's consumption in 2021 was approximately 76,389,561 MWh. 
However, part of Illinois' electricity already comes from either clean or renewable sources: 52.63$\%$ of it is generated by nuclear power plants \cite{Chicago_needs_58}. Another 12.27 $\%$ comes from wind energy, 1.50 $\%$ from solar, and 0.06$\%$ from hydroelectric plants \cite{Energy_Distribution}. Since there is no specific information on Chicago, we use these values. This means that we only need to cover for the 33.54$\%$ of the 76,389,561 MWh found above, which is 25,621,059 MWh. This is the value we use for our energy needs constraint - the minimum that needs to be produced in a year. This constraint appears in Inequality (\ref{energy_demand_4}) in the linear program in subsection \ref{section_4_linear_program_model}.

\subsubsection{Emissions} \label{section_4_constraints_emissions}

\quad Although both solar and wind systems generate no emissions during the production of energy, there is a quantifiable environmental impact from the production, transportation, and recycling of the materials used in plants. 

For wind plants, the production and transportation of the components, the reconditioning and renewal of the components, and disposal of the material are considered when quantifying the emissions of CO$_2$ \cite{WindEmissions_1}. On average, 4970g of CO$_2$ are emitted for every MWh generated by an onshore wind turbine, which is consistent with other studies \cite{WindEmissions_2}. 

For solar energy systems, silicon solar panels generate 45,000g CO$_2$ for every MWh when taking into account its production, maintenance and transportation \cite{SolarEmmissions_1}.

To set an upper bound for emissions, we estimate current emissions generated by non-renewable and non-clean sources. Of the 76,389,561 MWh consumed by Chicago, 20.99$\%$ comes from coal, 12.31$\%$ comes from natural gas, 0.21$\%$ - from biomass, and 0.04$\%$ - from oil \cite{Energy_Distribution}. The emissions resulting from these sources can be found in Table \ref{tab:current_emissions}.

{
\begin{table}[ht]
\centering
\begin{tabular}{ |p{2.5cm}||p{2cm}|p{3.5cm}|p{2cm}|  }
 \hline
 & Emissions per MWh (g CO$_2$ per MWh))&Energy produces a year in Chicago (MWh)&Total emissions (g CO$_2$)\\
 \hline
 Coal    & 820,000  & 16,034,169 (20.99$\%$) & 13.15 $\times$ 10$^{12}$\\

 Natural gas &   490,000  & 9,403,555 (12.31$\%$) & 4.61 $\times$ 10$^{12}$\\

 Oil (petroleum) &   1,106,765  & 30,556 (0.21$\%$) & 3.38 $\times$ 10$^{10}$\\

 Biomass &  230,000 & 160,418 (0.04$\%$) & 3.69 $\times$ 10$^{10}$\\
\hline
 Total &  &  & 17.83 $\times$ 10$^{12}$ \\
 \hline
\end{tabular}
\caption{Emissions from conventional energy sources \protect\cite{EmissionsOtherSources}, \protect\cite{EmissionsOtherSources_2}}
\label{tab:current_emissions}
\end{table}
} 

We estimate that 17.83 $\times$ 10$^{12}$ g CO$_2$ are being emitted every year by non-clean sources of energy to partially power the city of Chicago. We set our goal to reduce emissions by 80$\%$ to 3.578 $\times$ 10$^{12}$ g CO$_2$. This constraint appears in Inequality (\ref{emissions_4}) in the linear program in subsection \ref{section_4_linear_program_model}.

\subsubsection{Budget} \label{section_4_constraints_budget}

\quad This constraint is a limit on the funds allocated for the construction of new plants, since we assume that there is no infrastructure that is already operating and can be used. Although the following coefficients depend can vary, we chose the national average, shown in Table \ref{tab:cost_builduing}, which are obtained by dividing the average cost of a plant by its lifetime output in MWh.

{
\begin{table}[ht]
\centering
\begin{tabular}{|p{4cm}|p{5cm}|}
\hline
Levelized capital cost ($\$$ per MWh) & Source of energy \\ \hline
     27.45 &  wind (onshore)\\ 
      39.12  &  solar (hybrid) \\ \hline
\end{tabular}
\caption{Levelized capital cost per MWh \protect\cite{Cost_OM_1}, \protect\cite{Cost_OM_2}}
\label{tab:cost_builduing}
\end{table}
}

This constraint is shown in Inequality (\ref{budget_4}) in \ref{section_4_linear_program_model}.

\subsubsection{Space} \label{section_4_constraints_space}

\quad There is a limited amount of space available for new plants and infrastructures, and the Inequalities (\ref{space_4_1}, \ref{space_4_2}) reflect these limitations. We separate this space constraint into two inequalities. Solar systems have an advantage: they can be placed on the rooftops of existing building. Wind systems, on the other hand, need free space where plants can be built. 

\begin{figure}[ht]
    \centering
    \includegraphics[width=\textwidth]{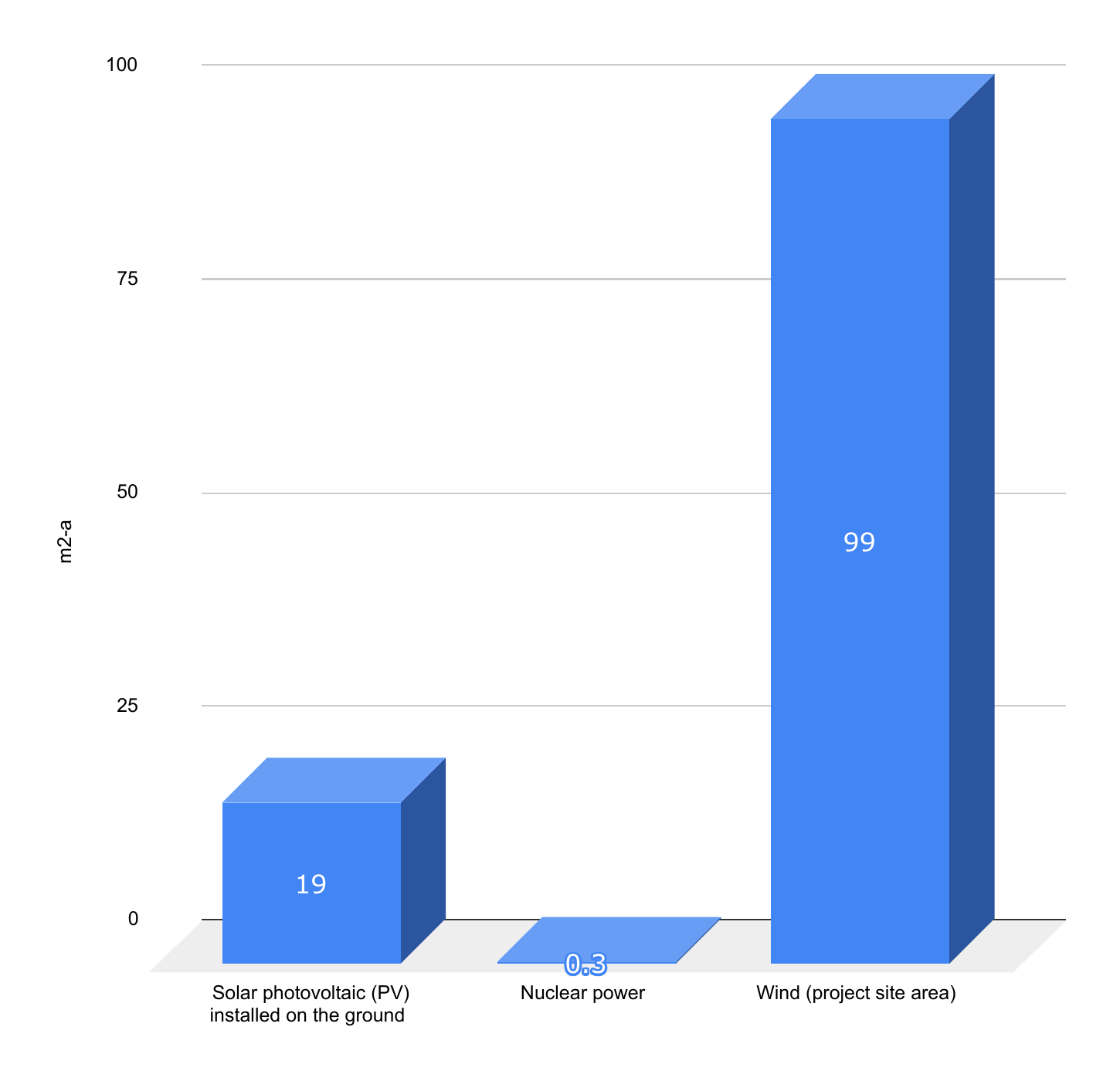}
    \caption{Average land-use in m$^2$a per every MWh produced  \protect\cite{space_UNECE_old}, \protect\cite{space_UNECE_new}, \protect\cite{Hannah_report}, \protect\cite{Hannah_calculations}}
    \label{fig:space-occupied}
\end{figure}

An estimation on land-use by different energy sources is shown in Figure \ref{fig:space-occupied}.

For solar panels, these occupy 19m$^2$ (204.5 ft$^2$) per MWh. The project site are for every MWh of wind energy is 99 m$^2$ (1065.6 ft$^2$) per MWh \cite{Hannah_calculations}.

For wind turbines, it is more difficult to approximate how much land can be used for new infrastructure. We estimate it as follows.

In 2010, the U.S. Census Bureau reported that 47$\%$ of total land in the United States was unoccupied \cite{Unoccupied_are_US_other}, \cite{Unoccupied_are_US_Wikipedia}. We apply this to the area of Illinois (1,614,570,000,000 ft$^2$), and get 758,847,900,000 ft$^2$. We choose to dedicate 1/15 of that ``unoccupied" area to wind power plants, which is 50,589,860,000 ft$^2$. This is our upper limit for space allocated for wind power plants, and thus we have the following constraint:

\begin{equation}
    1065.6 \frac{\text{ft$^2$}}{\text{MWh}}\times x_1 \leq 50,589,860,000 \text{ ft$^2$} \\
    \label{space_4_1_5}
\end{equation}

Dividing both sides by 1065.6 ft$^2$, we get an upper bound for wind energy production, which is the inequality (\ref{space_4_1}) in the model \ref{section_4_linear_program_model}.

For solar energy, we set a constraint that determines how much area of the existing buildings can be used to install solar panels. In this section we assume that solar panels can only be installed on the roofs. The building footprint in the city is approximately 70,532,107 ft$^2$ \cite{Buildings_area}, \cite{Buildings_area_not_accurate}. This is our upper limit for solar panels. 
\begin{equation}
    204.5 \frac{\text{ft$^2$}}{\text{MWh}}\times x_2 \leq 7.053 \times 10^7 \text{ ft$^2$} \\
\end{equation}

Dividing both sides by 204.5 ft$^2$, we get the upper bound for solar energy production, which is the inequality (\ref{space_4_2}) in 
\ref{section_4_linear_program_model}.

\subsection{Linear program} \label{section_4_linear_program_model}

\quad Minimize 
\begin{equation}
    C(x_1, x_2) = 37.80 x_1 + 58.62 x_2
\end{equation} 

\quad Subject to:
\begin{equation}
    x_1 + x_2 \geq 25,621,059  \text{MWh}
\label{energy_demand_4}
\end{equation}
\begin{equation}
    4970 x_1 + 45,000 x_2 \leq 3.578 \times 10^{12} \text{ g CO$_2$}
\label{emissions_4}
\end{equation}
\begin{equation}
     27.45 x_1 + 39.12 x_2 \leq  \$ 2 \times 10^9 
\label{budget_4}
\end{equation}
\begin{equation}
    x_1 \leq 47,475,469 \text{MWh}
\label{space_4_1}
\end{equation}
\begin{equation}
    x_2 \leq 344,900 \text{MWh}
\label{space_4_2}
\end{equation}
\begin{equation}
    0 \leq x_1 , x_2
\end{equation}\\

\subsection{Results} \label{section_4_results}
Using only solar and wind, the minimum of the objective function is obtained when exclusively wind energy is used. This occurred because wind turbines release fewer CO2 emissions during production than solar panels, cost less at comparable energy ratings, and operate at a higher efficiency. The results are presented in Table \ref{tab:linear_program_results_model_4}: 

{
\begin{table}[ht]
\centering
\begin{tabular}{ |p{3cm}||p{3cm}|p{1.5cm}|p{3cm}| }
\hline
& Wind energy&Solar energy&Total\\
\hline
Production (MWh)    & 25,621,059  & 0 & 25,621,059\\
\hline
Space occupied (ft$^2$) &    191,133,100  & 0 &  191,133,100\\
\hline
Emissions (g CO$_2$) &   127,336,663,230  & 0 & 127,336,663,230\\
\hline
Cost ($\$$) &  703,298,069  & 0 & 703,298,069 \\
\hline
Objective function values ($\$$) & 968,476,030 & 0 & 968,476,030 \\
\hline
\end{tabular}
\caption{Results from linear program \protect\ref{section_4_linear_program_model}}
\label{tab:linear_program_results_model_4}
\end{table}

}

This model's results are very important because:
\begin{itemize}
    \item It sets the minimum budget for the next three section: since wind is the cheapest of all energy sources to build and maintain, it means that if we try any other combination (but geothermal), it will be more expensive to build the necessary infrastructure. Below this value of 703,298,069 $\$$, the solution becomes unfeasible. 
    \item It also sets the minimum for emissions, in this case for the entire project. We will see later on that geothermal and nuclear energies have higher emissions per MWh, so this is the least amount we will be getting. Again, this tells us about the necessary conditions for feasibility of the solution.
\end{itemize}

\section{Energy Production and Demand throughout the day} \label{section_5_model_with_three_energy_constraints}

Since the production and demand depend on time of day, we split our energy needs constraint into three constraints which incorporate data from different periods of the day. In this section, we develop these new constraints and study the resulting model. 

\subsection{Demand and production throughout the day} 

\subsubsection{Production} \label{section_5_production_expl}

\begin{figure}[ht]
    \centering
    \includegraphics[width=0.95\textwidth]{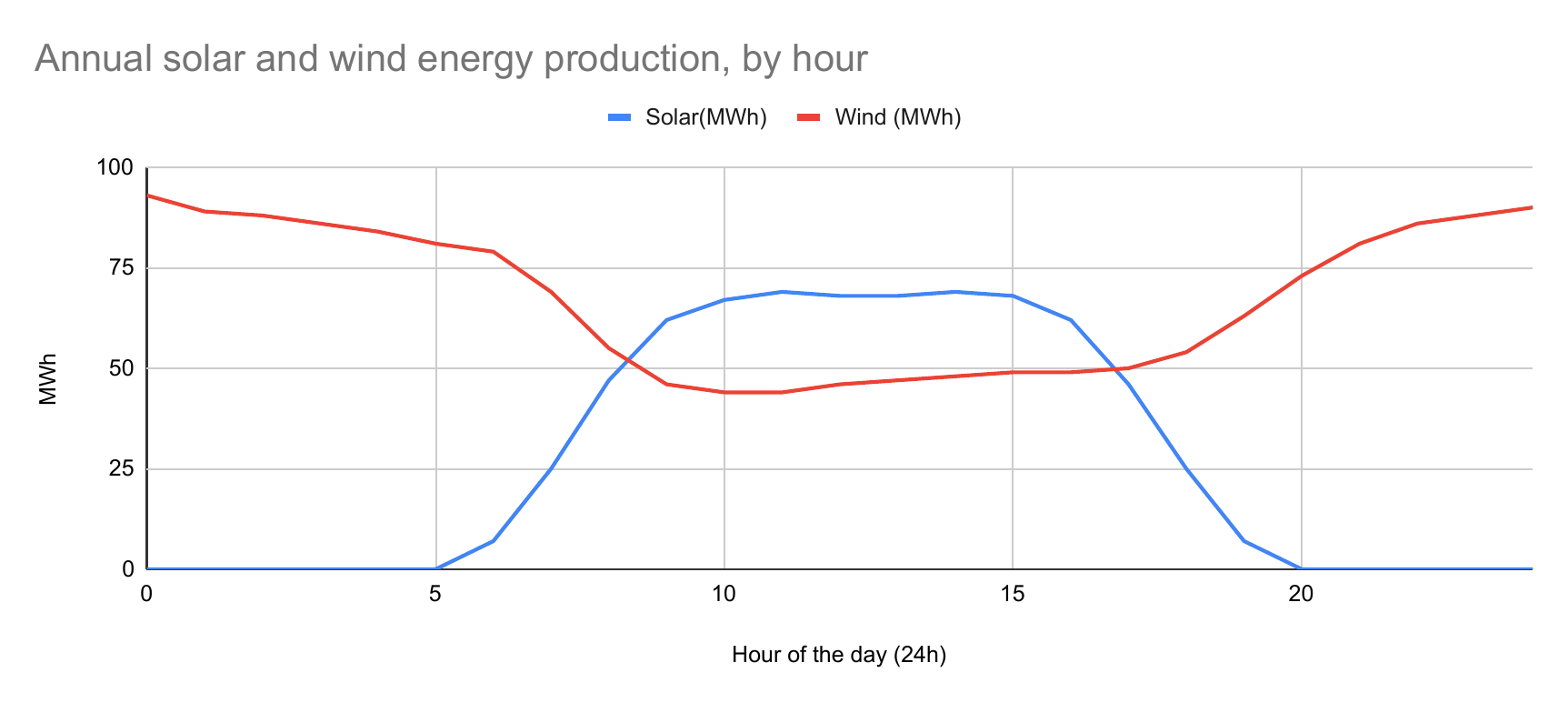}
    \caption{Energy Production by time of day: Solar (red) and Wind (blue). Solar only produces energy during the daytime, and Wind produces energy throughout the day with less production during daylight hours \protect\cite{Apex-Energy-Production}.}
    \label{fig:Apex-hourly-production}
\end{figure}

The outputs of wind and solar energy systems depend on the time of day. We split our single energy needs constraint in \ref{section_4_linear_program_model}, Inequality (\ref{energy_demand_4}), into three (\ref{energy_early_morning_5}, \ref{energy_daytime_5}, \ref{energy_evening_5}), one for each part time period of the day, shown in model \ref{section_5_linear_model_program}. The three periods of the day of choice were based on the output of solar systems, which only generate energy during the time of solar radiation, and this is approximately between 7am and 7pm (referred to as ``daytime hours" in this paper) \cite{Apex-Energy-Production}. The other two periods of time are the hours between 7pm and midnight (``evening hours"), and between midnight and 7am (``early morning hours"). Figure \ref{fig:Apex-hourly-production} shows the average production for both wind and solar energies throughout the day.

Based on this graph, the average production of both wind and solar systems are shown in Table \ref{tab:production_time_day}.

{
\begin{table}[ht]
\centering
\begin{tabular}{ |p{4cm}||p{3cm}|p{3cm}|  }
 \hline
 & Wind energy&Solar energy\\
 \hline
 Average production 12am-7am (MWh) and $\%$    & 600 MWh /37.69$\%$    &7 MWh /1.01$\%$\\
 \hline
 Average production 7am-7pm (MWh) and $\%$&   601 MWh /37.75$\%$  & 676 MWh /97.97$\%$\\
 \hline
 Average $\%$ of the daily demand &391 MWh /24.56$\%$ & 7 MWh /1.01$\%$\\

 \hline
\end{tabular}
\caption{Production based on the time of the day \cite{Apex-Energy-Production}}
\label{tab:production_time_day}
\end{table}

}

The percentages of wind and solar energy by time of day are used as coefficients in inequalities (\ref{energy_early_morning_5}, \ref{energy_daytime_5}, \ref{energy_evening_5}) in model \ref{section_5_linear_model_program}.

\subsubsection{Demand} \label{section_5_demand_expl}

Since our production is divided into early morning, daytime, and evening hours, we separate the data on demand similarly.

\begin{figure}[ht]
    \centering
    \includegraphics[width=\textwidth]{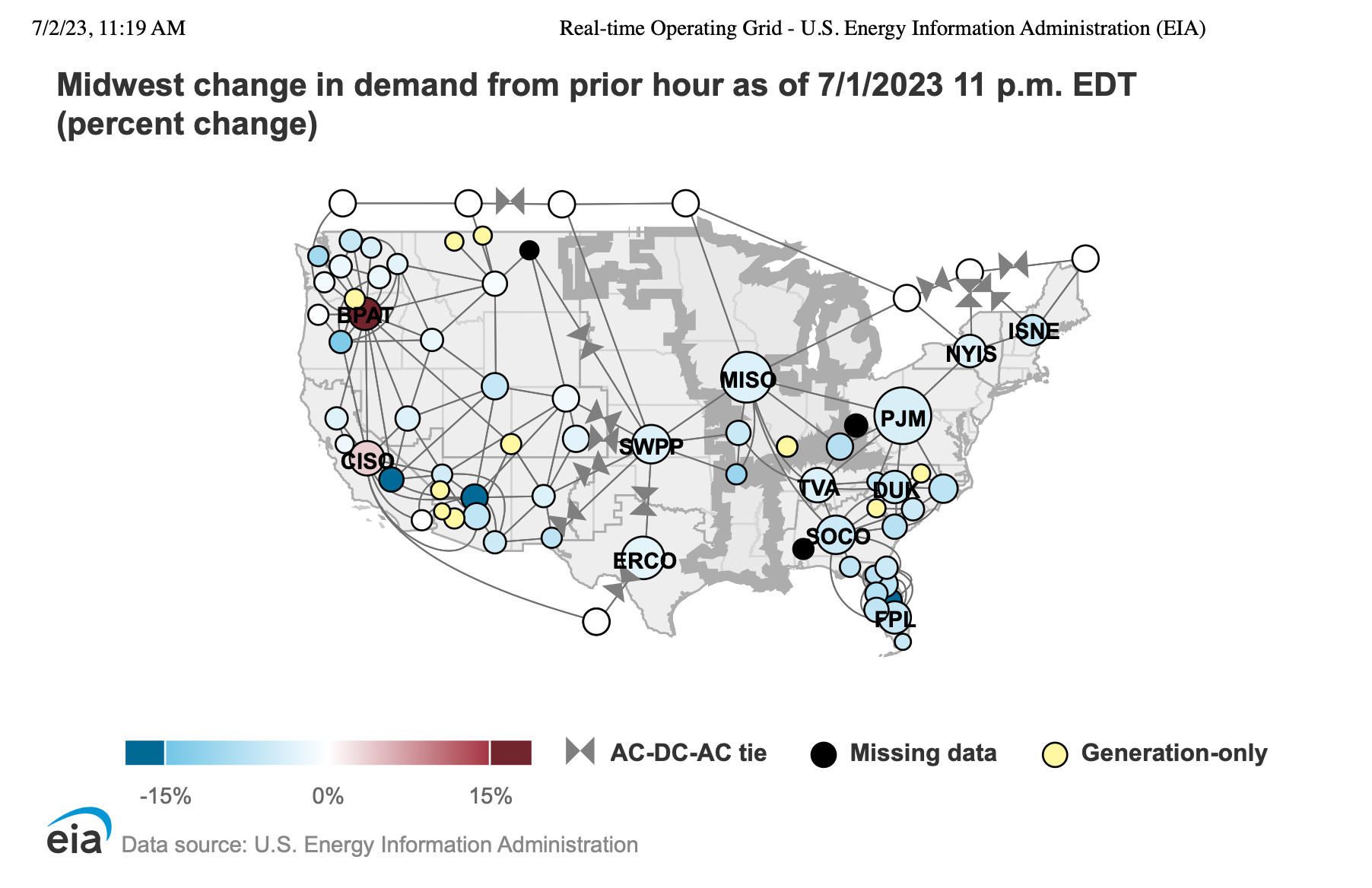}
    \caption{The data used to separate consumption on demand was obtained from the energy usage for the area highlighted in grey, the Midwest area \protect\cite{EnergyDemand}.}
    \label{fig:Energy_demand}
\end{figure}

In order to determine how much energy is used throughout the day and how that changes, we looked at the data from \cite{EnergyDemand} between July 2022 and June 2023, which provides demand on every hour of every day. Since no hourly energy usage of Chicago is available, we used data on the Midwest area, shown in Figure \ref{fig:Energy_demand}. We found an average over the early morning, daytime, and evening hours over the course of the year. We also calculated an average over these time periods for each month, and identified the largest percentages for each time period. The results are shown in the Table \ref{tab:demand_time_day}.

{
\begin{table}[ht]
\centering
\begin{tabular}{ |p{3cm}||p{2.5cm}|p{2.5cm}|p{2.5cm}|  }
 \hline
 & Early morning (12am-7am)&Daytime (7am-7pm)&Evening (7pm-12am)\\
 \hline
 Average $\%$ of the daily demand & 26.38$\%$ & 51.21$\%$&  22.41$\%$ \\
 \hline
 Highest $\%$ of the daily demand    &27.59$\%$ & 51.49$\%$&  23.44$\%$\\
 \hline
\end{tabular}
\caption{Demand based on time of the day \protect\cite{EnergyDemand}}
\label{tab:demand_time_day}
\end{table}

}

What we can see from this table is the following: 
\begin{enumerate}
\item 
Average $\%$ of the daily demand: this is simply telling us what percentages of energy is demanded during each of these day section of the day. On average, 26.38$\%$ is consumed in the early morning hours, 51.21$\%$ during the daytime hours, and 22.41$\%$ during the evening hours (for reference, 100$\%$ is the consumption during 24h).

\item 
Highest $\%$ of the daily demand: the previous row has the demand percentages obtained by averaging out the data over the entire year. However, there are some months when there is more need during the day than is shown, and similarly this can happen for the other two periods of the day. To make sure enough energy is produced every month for every section, we take the highest percentage. Although this implies that we will be overproducing, it also guarantees enough energy for every month of the year. 

\end{enumerate}

From this table, the important results are the percentages, which will help us to identify the lower bounds for (\ref{energy_early_morning_5}, \ref{energy_daytime_5}, \ref{energy_evening_5}). These constraints are formalized in section \ref{section_5_new_constraints_demand}.

\subsection{New constraints} \label{section_5_new_constraints_demand}

Based on the percentages we obtained in \ref{section_5_production_expl} and \ref{section_5_demand_expl}, the new energy need constraints (\ref{energy_early_morning_5}, \ref{energy_daytime_5}, \ref{energy_evening_5}) are: \\

\begin{itemize}
    \item Early morning needs (12am-7am): (25,621,059 MWh) $\times$ 0.2759 = 7.069 $\times $ 10$^6$ MWh;
    \item Day needs (7am-7pm): 25,621,059 MWh $\times$ 0.5149 = 13.192 $\times $ 10$^6$ MWh;
    \item Day needs (7pm-12am): 25,621,059 MWh $\times$ 0.2344 = 6.006 $\times $ 10$^6$ MWh.
\end{itemize}

Here is the updated model with three constraints for energy needs replacing the one constraint (\ref{energy_demand_4}) in \ref{section_4_linear_program_model}:\\

\subsection{Updated model} \label{section_5_linear_model_program}
\quad Minimize 
\begin{equation}
    C(x_1, x_2) = 37.80 x_1 + 58.62 x_2
\end{equation} 

\quad Subject to:
\begin{equation}
    0.3760 \times x_1 + 0.01 \times x_2 \geq 7.069 \times 10^6  \text{MWh}
\label{energy_early_morning_5}
\end{equation}
\begin{equation}
    0.3775 \times x_1  + 0.9797 \times x_2 \geq 13.192 \times 10^6 \text{MWh}
\label{energy_daytime_5}
\end{equation}
\begin{equation}
    0.2456 \times x_1 + 0.01 \times x_2  \geq 6.006 \times 10^6 \text{MWh}
\label{energy_evening_5}
\end{equation}
\begin{equation}
    4970 x_1 + 45,000 x_2 \leq 3.578 \times 10^12  \text{ g CO$_2$}
\label{emissions}
\end{equation}
\begin{equation}
    27.45 x_1 + 39.12 x_2 \leq \$ 2 \times  10^9 
\label{budget_5}
\end{equation}
\begin{equation}
    x_1 \leq 47,475,469 \text{MWh}
\label{space_5_1}
\end{equation}
\begin{equation}
    x_2 \leq 344,900 \text{MWh}
\label{space_5_2}
\end{equation}
\begin{equation}
    0 \leq x_1 , x_2
\end{equation}\\

\subsection{Results} \label{section_5_results}

{
\begin{table}[ht]
\centering
\begin{tabular}{ |p{3cm}||p{2.5cm}|p{2.5cm}|p{3cm}|  }
 \hline

 & Wind energy&Solar energy&Total\\
 \hline
 Production 12am-7am (MWh)    &  12,959,826  & 3,449& 12,963,275\\
 \hline
 Production 7am-7pm (MWh)    & 12,854,101  & 337,898 & 13,192,000\\
 \hline
 Production 7am-7pm (MWh)    & 8,185,153   & 3,449 &  8,188,602\\
 \hline
 Production total (MWh)    & 34,104,806   &344,900&  34,449,706\\
 \hline
 Emissions (g CO$_2$) &   169,500,885,820  & 15,520,500,000 &  185,021,385,820\\
 \hline
 Budget ($\$$) &  936,176,924 & 13,492,488 &  949,669,412\\
 \hline
 Objective function values ($\$$) & 1,289,161,666 & 20,218,038 &  1,309,379,704 \\
 \hline

\end{tabular}
\caption{Results of the linear program \protect\ref{section_5_linear_model_program}.}
\label{tab:results_section_daytime_separate}
\end{table}
} 

First we look at other constraints and the values we got: 
\begin{itemize}
    \item Emissions: we can see that our current emissions are less than 5$\%$ of our upper limit and 1.04$\%$ of the emissions that would have been produced if we chose non-renewable energies. 
    \item Budget: the total price to construct enough plants to get our energy needs met is 949,669,412 $\$$ a year. 
\end{itemize}

Next, we look at the energy production. As expected, we are overproducing in general: our needs for Chicago are 25,621,059 MWh, while our total production given by the linear program is 34,449,706 MWh. This is partially due to our choice of picking the highest percentages of demand for each section of the day. Notice that we are producing the maximum amount of solar energy that is allowed in Inequality (\ref{space_5_2}). 
In section \ref{section_6_space_solar}, we assume that solar panels can be installed on the ground as well, and study the resulting linear program.

\section{Allocating more space for solar systems}  \label{section_6_space_solar}

In this section, we allow for solar energy system installation on the ground. Remember that solar systems use 204.5ft$^2$ for every MWh they produce when installed on the ground \cite{space_UNECE_new}. Now we add it to our space constraint (\ref{space_4_1_5}):

\begin{equation}
    1065.6 x_1 + 204.5 x_2 \leq 50,589,860,000 ft^2
    \label{almost_done_space_constraint}
\end{equation}

The Inequality (\ref{almost_done_space_constraint}) is going to substitute Inequalities (\ref{space_5_1}) and (\ref{space_5_2}) from \ref{section_5_linear_model_program}, because solar panels are now sharing space with wind plants on the ground. 
We also need to account for the 344,900 MWh of energy from solar panels installed on top of buildings, found in section \ref{section_4_constraints_space}. Adding it to the new space constraint Inequality (\ref{almost_done_space_constraint}), we get:

\begin{equation}
    1065.6 x_1 + 204.5 (x_2 - 344,900 MWh) \leq 50,589,860,000 ft^2
\end{equation}

\subsection{New model with space constraint on solar energy production} \label{section_6_linear_model_program}

\quad Minimize 
\begin{equation}
    C(x_1, x_2) = 37.80 x_1 + 58.62 x_2
    \label{old_objective_function}
\end{equation} 
\quad Subject to:
\begin{equation}
    0.3760 \times x_1 + 0.01 \times x_2 \geq 7.069 \times 10^6 \text{MWh}
\label{energy_early_morning_6}
\end{equation}
\begin{equation}
    0.3775 \times x_1 + 0.9797 \times x_2 \geq 13.192 \times 10^6 \text{MWh}
\label{energy_day_6}
\end{equation}
\begin{equation}
    0.2456 \times x_1 + 0.01 \times x_2 \geq 6.006 \times 10^6 \text{MWh}
\label{energy_evening_hours_6}
\end{equation}
\begin{equation}
    4970 x_1 + 45,000 x_2 \leq 16,325 \times 10^9 \text{ g CO$_2$}
\label{emissions_6}
\end{equation}
\begin{equation}
    27.45 x_1 + 39.12 x_2 \leq \$ 2 \times  10^9 
\label{budget_6}
\end{equation}
\begin{equation}
    1065.6 x_1 + 204.5 (x_2 - 2,190,438 \text{MWh}) \leq 50,589,860,000 \text{ ft$^2$}
\label{space_6}
\end{equation}
\begin{equation}
    0 \leq x_1 , x_2
\end{equation}\\

\subsection{Results} \label{section_6_results}

{
\begin{table}[ht]
\centering
\begin{tabular}{ |p{3cm}||p{2.5cm}|p{2.5cm}|p{3cm}|  }
 \hline

 & Wind energy&Solar energy&Total\\
 \hline
 Production 12am-7am (MWh)    & 9,447,742   & 39,005& 9,486,747\\
 \hline
 Production 7am-7pm (MWh)    & 9,370,668  & 3,821,331& 13,192,000\\
 \hline
 Production 7pm-12am (MWh)    & 5,966,994   & 39,053 &  6,006,000\\
 \hline
 Production total (MWh)    & 24,862,479   &3,900,512&  28,762,991\\
 \hline
 Space occupied (ft$^2$) &   26,493,457,622 & 797,309,844  &  27,290,767,467\\
 \hline
 Emissions (g CO$_2$) &   123,566,520,630  & 175,523,049,000 &  299,089,569,630\\
 \hline
 Budget ($\$$) &   682,475,048 & 152,588,037 &  835,063,085\\
 \hline
 Objective function values ($\$$) & 939,801,706 & 228,648,025 &  1,168,449,731 \\
 \hline

\end{tabular}
\caption{Results of the linear program \protect\ref{section_6_linear_model_program}}
\label{tab:my_label}
\end{table}
} 

Notice that solar energy production increased considerably, and the value of the objective function decreased by over $\$$150 million.

However, we are still overproducing, in particular during the early morning hours. Allocating more space for solar panels made this model more successful, but did not completely solve the issue of overproduction. Two possible solution that we consider in this project are: considering storing the excess energy (which we do no explore due to the scope of this project), or adding another source of energy, such as nuclear (section \ref{section_7_nuclear}) or geothermal (section \ref{section_8_geothermal}).

\section{Introducing nuclear energy into the model} \label{section_7_nuclear}

In this section, we explore nuclear power as a solution to the issue of  overproduction. 

\subsection{Constraints for nuclear power} \label{section_7_constraints_nuclear_model}

\subsubsection{Energy production constraint}
If we choose to use nuclear power, the smallest unit of system we can construct is one reactor. The smallest nuclear reactors (called Small Modular Reactors, SMRs) have the capacity of 300 MW \cite{Smallest_reactor}, although an average reactor would be of 1 GW \cite{Average_capacity}. A plant with an SMR produces 2,628,000 MWh a year working at full capacity, and we will set this as our lower bound for nuclear energy production, as shown in Inequality (\ref{minimum_nuclear}) in \ref{section_7_linear_program_model}.

\subsubsection{Time of the day constraint} \label{section_7_production_nuclear}

An important advantage of nuclear power is that its output distribution (the amount of energy it delivers) is constant, and can also be flexible \cite{NucEnergyOutput}. In our model, we assume that the energy production of a power plant is constant throughout the day, meaning that 29$\%$ is produced during the early morning hours (\ref{early_morning_7}), 50$\%$ during the daytime hours (\ref{daytime_7}), and 21$\%$ during the evening hours (\ref{evening_7}) (based on the fact that early morning hours (7 hours) are 29$\%$ of the day, daytime (12 hours) are 50$\%$, etc). These coefficients are reflected in Inequalities (\ref{early_morning_7}, \ref{daytime_7}, \ref{evening_7}) in \ref{section_7_linear_program_model}

\subsubsection{Emission constraint} \label{section_7_emissions_nuclear_model}
Although nuclear energy is considered a clean energy source because there are no direct emission during nuclear fission, it is not a renewable source because of its use of uranium. So, there are emissions produced during the construction and operation of the plant itself, and the mining of uranium.

Carbon dioxide is emitted during the mining of uranium. On average, 34,000 gCO$_2$ are estimated to be emitted for every MWh produced using uranium, and this number is predicted to increase as the amount of uranium left decreases \cite{UraniumTotal}. Also, it is important to mention that for every kWh produced from a nuclear fission, 0.1-0.3 kWh are used \cite{NucNeeds}. We will use 0.2 kWh in our calculations. This means that to obtain 1 MWh of electricity, we will emit 40,800 gCO$_2$. 

Constructing a plant and maintaining it is associated with emissions of GHG as well. An average of 8200 gCO$_2$ are emitted for every MWh that a plant produces \cite{NucConstEmissions}.

To get the final coefficient for the emissions constraint, we add the two emission values: for uranium and for the construction/maintenance. These add up to 49,000 g CO$_2$/MWh, which is reflected in Inequality (\ref{emissions_7}) in \ref{section_7_linear_program_model}.

\subsubsection{Budget constraint} \label{section_7_budget_nuclear}
The cost for every additional MWh produced by a nuclear reactor is approximately 70.8$\$$ \cite{Cost_OM_2}. This is shown in Inequality (\ref{budget_7}) in \ref{section_7_linear_program_model}.

\subsubsection{Space constraint} \label{space_nuclear}
As reported in \cite{space_UNECE_new} and shown in figure \ref{fig:space-occupied}, 0.3 m$^2$ are needed for every MWh, which is 3.23 ft$^2$, shown in Figure \ref{fig:space-occupied}. This coefficient is represented in Inequality (\ref{space_7}) in \ref{section_7_linear_program_model}

\subsection{Model} \label{section_7_linear_program_model}
The adjusted model in this section is composed of our previous model of wind energy (x$_1$) and solar energy (x$_2$) from Section \ref{section_5_model_with_three_energy_constraints} with the addition of nuclear energy (x$_3$). 

\quad Minimize 
\begin{equation}
    C(x_1, x_2, x_3) = 37.80 x_1 + 58.62 x_2 + 96.2 x_3
\end{equation} 

\quad Subject to:
\begin{equation}
    0.3769 x_1 + 0.01 x_2 + 0.29 x_3 \geq 7.069 \times 10^6  \text{MWh}
    \label{early_morning_7}
\end{equation}
\begin{equation}
    0.3775 x_1 + 0.9797 x_2 + 0.5 x_3 \geq  13.192 \times 10^6 \text{MWh}
    \label{daytime_7}
\end{equation}
\begin{equation}
    0.2456 x_1 + 0.01 x_2 + 0.21 x_3 \geq  6.006 \times 10^6 \text{MWh}
    \label{evening_7}
\end{equation}
\begin{equation}
    4970 x_1 + 45,000 x_2 + 49,000 x_3 \leq 163,325 \times 10^9 \text{ g CO$_2$}
    \label{emissions_7}
\end{equation}
\begin{equation}
    27.45 x_1 + 39.12 x_2 + 70.8 x_3 \leq \$ 2 \times  10^9
    \label{budget_7}
\end{equation}
\begin{equation}
    1065.6 x_1 + 204.5 (x_2 - 10,279,088) + 3.23 x_3 \leq 50,589,860,000 \text{ ft$^2$}
    \label{space_7}
\end{equation}
\begin{equation}
    0 \leq x_1, x_2
\end{equation}
\begin{equation}
    2,628,000 \leq x_3
    \label{minimum_nuclear}
\end{equation}
\\

\subsection{Results} \label{section_7_results}

First, we ran the linear program as it is shown in the subsection \ref{section_7_linear_program_model}. Given these constraints, the result that minimized the objective function was the same as in section \ref{section_6_results}. In other words, the program suggests to only use wind and solar systems to cover the energy needs. This is a somewhat expected result - nuclear power is more expensive than wind and energy to both maintain and construct. 

So, why would anyone choose nuclear power instead of wind or solar? If you look at the constraints, the answer becomes obvious: land use. Nuclear power plants occupy less space per MWh than solar or wind plants. In fact, if the space constraint is reduced to something lower than what we chose in section \ref{section_4_constraints_space}, then the linear program does suggest we use nuclear power. 

To demonstrate this, we reduced the upper limit on space to 205,898,600 ft$^2$,and obtain the following constraint on space:\\
\begin{equation}
    1065.6 x_1 + 204.5 (x_2 - 10,279,088) + 3.23 x_3 \leq 205,898,600 \text{ft$^2$}
    \label{section_7_space_new_nuclear}
\end{equation}

The results obtained with this new constraint on space are shown in Table \ref{tab:nuclear_results}.

{
\begin{table}[ht]
\centering
\begin{tabular}{ |p{2cm}||p{2cm}|p{2cm}|p{2cm}|p{2.5cm}|  }
 \hline

 & Wind energy&Solar energy&Nuclear energy&Total\\
 \hline
 Production 12am-7am (MWh)    & 9,812,073  & 21,904 & 762,120 & 10,596,098\\
 \hline
 Production 7am-7pm (MWh)    &   9,732,027 &  2,145,972 & 1,314,000 & 13,192,000\\
 \hline
 Production 7pm-12am (MWh)    & 6,197,099  & 21,904 & 551,880 & 6,770,883\\
 \hline
 Production total (MWh)    & 25,821,247  & 2,190,438 & 2,628,000  & \\
 \hline
 Space occupied (ft$^2$) &  192,626,502  &0  & 8,488,440 & 201,114,942\\
 \hline
 Emissions (g CO$_2$) &  128,331.6 $\times$ 10$^6$ &  98,569.71 $\times$ 10$^6$ & 128,772 $\times$ 10$^6$ & 355,673,307,590\\
 \hline
 Cost ($\$$) &  708,793,230   & 85,689,934 & 186,062,400 & 980,545,564\\
 \hline
 Objective function values ($\$$) & 976,3043,136 & 128,403,475 & 252,813,600 & 1,357,260,212\\
 \hline

\end{tabular}
\caption{Results of the linear program \protect\ref{section_7_linear_program_model} with adjusted space constraint (\protect\ref{section_7_space_new_nuclear}).}
\label{tab:nuclear_results}
\end{table}
}

Notice that we are very close to the upper limit for our space constraint. The amount of solar energy we are producing is what we can ``fit" on the buildings, but it does not occupy any land. Moreover, as we increase the upper limit for the space constraint (\ref{space_7}), the linear program suggests to increase the amount of MWh produced by wind systems or solar systems and reduce nuclear production. In other words, we only use nuclear power when we do not have enough space to build another system and run a different source of energy.

Furthermore, nuclear power does not solve the issue of overproduction. The daytime needs are met, but a large amount of energy is produced during the morning hours and the evening hours. We conclude that nuclear energy is not the best solution unless we lack space to install other systems. However, the advantage that nuclear energy has over the rest is that it does not depend on either the weather or natural resources, which means that it is the only clean energy source that can ensure a constant supply.

\section{Geothermal energy} \label{section_8_geothermal}
In this section, we consider adding geothermal energy to the linear program \ref{section_6_space_solar}. Since nuclear energy did not reduce the value of the objective function and increased some of the constraints, as shown in Table \ref{tab:nuclear_results}, it will not be included in this section.

\subsection{Constraint for geothermal energy} \label{section_8_constrints_geo}

\subsubsection{Constraint on emissions} \label{section_8_geothermal_emissions}

Although geothermal energy is considered a clean energy, the systems require raw materials, transportation, etc. Geothermal systems emit, as a median, 38,000 gCO$_2$ per every MWh produced \cite{EmissionsOtherSources}. This is reflected in Inequality (\ref{emissions_geothermal}) in \ref{section_8_linear_program_model}.

\subsubsection{Constraint on space} \label{section_8_geothermal_space}

Geothermal systems occupy a considerable amount of space \cite{Geothermal_plant_structure}.
Around 900 m$^2$ are required to produce every additional GWh  \cite{Geothermal_main_source} by a 56 MW geothermal flash plant. This  is equivalent to 9.6875 ft$^2$/MWh, which is added to the space constraint in Inequality (\ref{space_geothermal}) in \ref{section_8_linear_program_model}.

\subsubsection{Constraints based on the time of the day} \label{section_8_geothermal_energy_production_daytime}

Geothermal energy, just like nuclear energy, has a great advantage - power plants produce electricity consistently, running 24 hours a day, 7 days a week \cite{Geothermal_production_time_day}. This implies that a plant produces approximately 4.167$\%$ of its total daily output during every hour of the day. Multiplying this by the amount of hours in each of our day sections (12am-7am, 7am-7pm, 7pm-12am), we get the following coefficients:\\
\begin{itemize}
    \item Early morning hours: 7h $\times$ 4.167$\%$ = 29.16$\%$, or 0.2916;
    \item Day hours: 12h $\times$ 4.167$\%$ = 50$\%$, or 0.5000;
    \item Evening hours: 5h $\times$ 4.167$\%$ = 20.83$\%$, 0.2083.
\end{itemize}

These coefficients that reflect the energy production of a geothermal plant are added to the three constraints on energy production in \ref{section_6_space_solar}, and can be found as Inequalities (\ref{early_morning_geothermal}, \ref{daytime_geothermal}, \ref{evening_geothermal}) in \ref{section_8_linear_program_model}.

\subsection{Model} \label{section_8_linear_program_model}

Here is the new linear program that includes wind energy (x$_1$), solar energy (x$_2$), and geothermal energy (x$_4$), all in MWh. 

\quad Minimize 
\begin{equation}
    C(x_1, x_2, x_4) = 73.7 x_1 + 55.8 x_2 + 39.61 x_4
\end{equation} 

\quad Subject to:
\begin{equation}
    0.3769 x_1 + 0.01 x_2 + 0.2916 x_4 \geq 7.069 \times 10^6  \text{MWh}
    \label{early_morning_geothermal}
\end{equation}
\begin{equation}
    0.3775 x_1 + 0.9797 x_2 + 0.5 x_4 \geq  13.192 \times 10^6 \text{MWh}
    \label{daytime_geothermal}
\end{equation}
\begin{equation}
    0.2456 x_1 + 0.01 x_2 + 0.21 x_4 \geq  6.006 \times 10^6 \text{MWh}
    \label{evening_geothermal}
\end{equation}
\begin{equation}
    4970 x_1 + 45,000 x_2 + 38,000 x_4 \leq 3.578 \times 10^12 \text{ g CO$_2$}
    \label{emissions_geothermal}
\end{equation}
\begin{equation}
     27.45 x_1 + 39.12 x_2  + 21.8 x_4 \leq \$ 2 \times  10^9 
    \label{budget_geothermal}
\end{equation}
\begin{equation}
    1065.6 x_1 + 204.5 (x_2 - 10,279,088) + 9.6875 x_4 \leq 50,589,860,000 \text{ ft$^2$}
    \label{space_geothermal}
\end{equation}
\begin{equation}
    0 \leq x_1, x_2, x_4
\end{equation}\\

\subsection{Results} \label{section_8_results}

{
\begin{table}[ht]
\centering
\begin{tabular}{|p{2.5cm}||p{2.5cm}|p{1cm}|p{2.5cm}|p{2.5cm}| }
 \hline

 & Wind energy&Solar energy&Geothermal energy&Total\\
 \hline
 Production 12am-7am (MWh)    & 2,164,412  & 0 & 6,406,242 & 8,570,654\\
 \hline
 Production 7am-7pm (MWh)    &   2,146,755 &  0 & 11,045,245 & 13,192,000\\
 \hline
 Production 7pm-12am (MWh)    & 1,366,997  & 0 & 4,639,002& 6,006,000\\
 \hline
 Production total (MWh)    & 5,695,821  & 0 &  22,090,490  & \\
 \hline
 Space occupied (ft$^2$) &   42,490,827 &0  & 214,001,621 & 256,492,449\\
 \hline
 Emissions (g CO$_2$) &  28,308,232,358  &  0 & 839,438,620,000 & 867,746,852,358\\
 \hline
 Cost ($\$$) &  156,350,297   & 0 & 481,572,682 & 637,922,979\\
 \hline
 Objective function values ($\$$) & 215,302,048 & 0 & 971,981,560 & 1,090,306,357 \\
 \hline

\end{tabular} 
\caption{Results of the linear program including geothermal energy \protect\ref{section_8_linear_program_model}}
\label{tab:geothermal_results}
\end{table}} 

This is the lowest value of the objective function we have been able to achieve, and it is the result of geothermal energy being cheaper than solar. Notice that the emissions, although still within the limits, have increased, because we are using less wind energy, which is 10 times cleaner than geothermal. 

\section{Discussion}

Current energy production methods are unsustainable in the long run.
A good substitution for these methods are renewable energy, which has lower emissions and are sustainable over time, and clean energy, which produces fewer emissions as well. 
In this paper, we study minimization of total renewable and clean energy costs with constraints on emissions.
We treat total energy cost as a linear function and establish linear constraints from practical consideration.
The resulting linear programs can be solved using the simplex method. 

We constrain the problem by requiring that energy needs are met.
If meeting energy needs is the only concern, it is most cost effective to use natural gas; however, natural gas comes with high emissions. 
We explore alternative sources to reduce emissions while still meeting energy needs.


In this paper, we have studied a total of five models which address different concerns when using clean and renewable energies.
We have considered three models that combine wind and solar energy sources with constraints on emissions, budget, and space.  
From the simplest model  
\ref{section_4_linear_program_model}, in section \ref{section_4_model_using_solar_and_wind_no_space_constraint}, we have learned that wind plants are overall the cleanest source of renewable energy, as well as the cheapest to build, operate and maintain.
However, wind and solar plants' output is not constant throughout the day. Additionally, the demand on energy is not constant either. This fact warranted a division of the day (24h) in 3 slots (12am-7am, 7am-7pm, and 7pm-12am). The hours were divided based on the output of solar energy systems. After adjusting for demand and production rates during these three slots, we found from a second model in section \ref{section_5_model_with_three_energy_constraints} that the minimum value for the objective function is achieved as a result of a combination of solar and wind systems. In the third model in \ref{section_6_space_solar}, we allowed for land use of solar systems, which led to an increase of suggested amount of MWh produced by solar panels and, subsequently, a decrease in wind energy production. 
However, we were facing the problem of overproduction: wind plants were producing more energy than needed during some slots, and because wind energy cannot be easily stored for later use, this resulted in wasted MWh. To mitigate this issue, we considered two solutions: adding nuclear energy production and geothermal energy production. 
%
In section \ref{section_7_nuclear}, we found from a fourth model \ref{section_7_linear_program_model} that, although nuclear energy reduced the overproduction slightly, it was not financially optimal given that nuclear plants are more expensive to build and operate. Adding nuclear sources also increased emissions. We concluded that nuclear plants were not a good solution for overproduction unless the available space was drastically reduced. 
In section \ref{section_8_geothermal}, we studied a fifth model with added geothermal plants to the system. This was more successful than nuclear in mitigating overproduction, and it resulted in a reduced objective function value. Geothermal energy also eliminated the need for solar energy and reduced the wind energy production, although this came with a cost: slightly increased emissions.
Overall, a combination of geothermal energy and wind energy is the best financial decision. They can meet the emissions, budget, and space constraints while providing the necessary electricity levels and minimizing overproduction. 

The results of each model were obtained given an objective function that considers building, operation and maintenance costs. As mentioned in \ref{section_4_constraints_budget}, we assume that all the infrastructure needs to be built/purchased. But what about preexisting plants in Illinois available for use? We studied the case in which the objective function coefficients for each energy source excluded the building costs, keeping the maintenance and operation expenses. The change of the objective function did not yield a different balance of energy sources, which can be seen explicitly for the model from section~\ref{section_6_space_solar} in Appendix, section~\ref{change_objective_function_section}.

Another aspect of this study that was considered is: what if we chose to minimize emissions instead of costs? Since wind energy has the lowest emissions per MWh by far, the minimal emissions are obtained from using only wind energy,  Appendix \ref{minimizing_emissions_constraint}.

This study is not without caveats. There are some limitations associated with the energy sources of choice and the results we have obtained. All the energy sources considered, except for nuclear, depend on the natural resources and weather conditions of their building site. In particular, the maximum output of wind, solar, and geothermal systems varies significantly based on their location \cite{efficiency_renewables}. Given that the information available on the capacity factor \cite{Capacity_factor_dict} is limited, we used the average (across the U.S. for some and the world for other) values for our models. This approximation causes both over- and underestimation in our results, which should be addressed once more information on capacity factors for different regions are available. Currently available capacity factors for different regions should be used instead of the generalized ones when applying the model to other cities. For instance, Phoenix is expected to have higher solar energy output than Chicago, which may not be the case for wind energy. 

Another limitation of this study that deserves more attention is the energy needs constraints. Given the scope of this study, we chose to have three constraints for energy needs based on the varying energy outputs throughout the day. However, the outputs also vary depending on the month of the year. This is especially important for those cities with larger weather variability, like New York or Chicago. To achieve more accurate constraint bounds, it is wise to consider creating models based on the month or even week of the year. 

In this study, we limited the constraints to budget, emissions, land use, and energy needs. However, there are more factors that come into play when choosing a city's sources of energy. For instance, the number of jobs new plants would generate (maintenance, construction) and eliminate (due to lack of need of non-renewable sources), the time new plants would take to be construct, availability of materials for the infrastructure and systems, changes in levelized capital costs due to technology advances and tax incentives, among other factors. 

Our project aimed to investigate the feasibility of fully replacing traditional energy sources with sustainable and clean alternatives to meet the essential energy demands of a city like Chicago while maintaining a realistic cost. While we successfully identified optimal combinations of renewable and clean energy sources to minimize costs while fulfilling these baseline requirements, questions lingered about the practicality of these findings. It would be very costly to substitute all non-renewable energy with renewable energy, especially given growing energy demands, due to the substantial infrastructure development required. On the other hand, clean energy options, such as nuclear power, rely on specific finite natural resources like uranium.


In essence, relying solely on renewable and clean energy sources to keep pace with growing energy demands appears unattainable due to the extensive infrastructure and resource requirements. Conversely, non-renewable energy sources are becoming scarcer and costlier to extract as easily accessible natural resources dwindle. Our energy supply is struggling to match the rapid growth of our energy consumption. Neither traditional nor renewable sources can adequately meet our escalating demands. To simultaneously mitigate environmental impact and secure our energy future, we must reduce our consumption. Ideally, the future of our energy lies in a combination of renewable and non-renewable resources, but maintaining a realistic baseline energy requirement is essential to make this blend sustainable.

\clearpage
\bibliographystyle{abbrv}
\bibliography{energybib}

\newpage
\appendix

\section{Change in the objective function}  \label{change_objective_function_section}

\subsection{Model} \label{linear_model_3_1}

Before we move on to solving the issue of overproduction, we would like to talk about the objective function we are using. The objective function is based two assumptions: 

\begin{itemize}
    \item There is no available infrastructure (wind plants, solar systems, etc) that can be used to cover Chicago's needs.
    \item There is no funding to build the plants, and that we will pay for the ``construction of every additional MWh" as it gets produced.
\end{itemize}
Hence, our budget constraint served as an upper limit to how much infrastructure we can afford to ``pay off". However, we can consider the case when we have a budget large enough to pay for all the infrastructure. In this case, we would be looking to minimize the operation and maintenance costs only. After finding the fixed O$\&$M cost values in \cite{Cost_OM_1} and \cite{Cost_OM_2}, we consider an alternative objective function. The new coefficients are shown in Table \ref{tab:change_objective_fuinction}.

{
\begin{table}[ht]
\centering
\begin{tabular}{|p{4cm}|p{5cm}|}
\hline
Levelized O$\&$M and transmission costs & kind of energy \\ \hline 
     10.35 &  wind (onshore)\\ 
      19.51  &  solar (hybrid) \\ \hline
\end{tabular}
\caption{Costs of operation and maintenance.}
\label{tab:change_objective_fuinction}
\end{table}
}

and the new objective function is:\\

\begin{equation}
    C(x_1, x_2) = 10.35 x_1 + 19.51 x_2
    \label{new_objective_function}
\end{equation}

Now we determine whether this changes the output of x$_1$ (wind energy production a year in MWh) and x$_2$ (solar production a year in MWh). The constraints will be the same as in the previous model: (\ref{energy_early_morning_6}, \ref{energy_day_6}, \ref{energy_evening_hours_6}, \ref{emissions_6}, \ref{budget_6}, \ref{space_6}).

\subsection{Results}

\begin{figure}[ht]
    \centering
    \includegraphics[width=\textwidth]{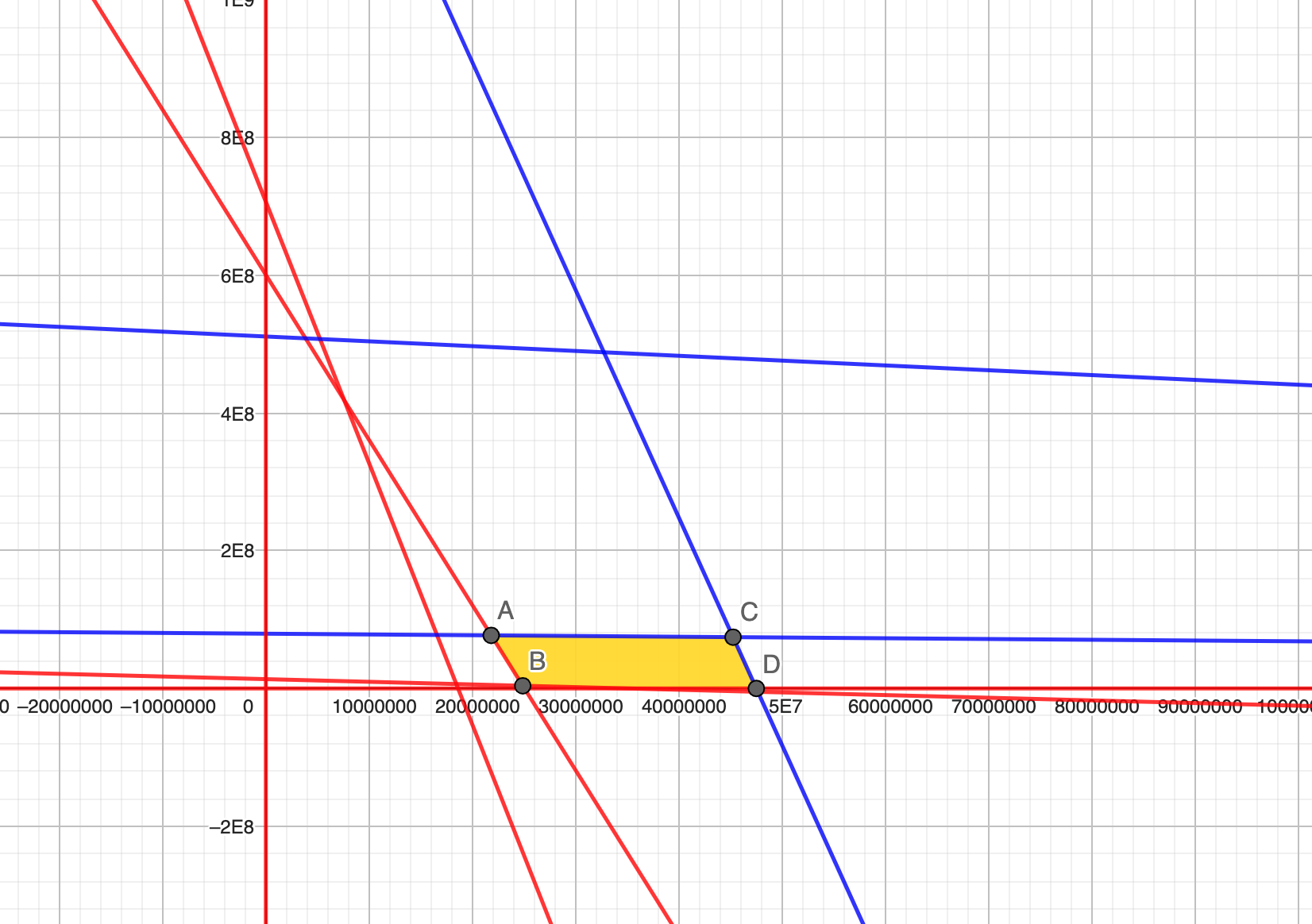}
    \caption{Constraints on the linear program from section \protect\ref{linear_model_3_1}. Red section represent the lower bounds (energy needs), and blue section are the upper bounds (emissions, budget, and space). The yellow polygon represents all possible solutions, with A, B, C, and D being the corner points. The axis x:y are scaled to 1:15.}
    \label{fig:constraints_colors}
\end{figure}

The results of this program turned out to be the exact same as of the model in section 6.2. The total production of wind energy was 24,862,479 MWh, and solar energy was 3,900,512 MWh. Consequently, all the constraint values were the same. 
Why didn't the result change? Visualizing the problem will help us understand it mathematically. 
When the feasible region of a linear program is bounded, it is a convex polygon and the minimum of the objective function occurs at a vertex (James K Strayer text).  
Refer to Figure \ref{fig:constraints_colors}, where horizontal axis represent the MWh of wind energy, and the vertical axis - MWh of solar energy.

The four corner points from Figure \ref{fig:constraints_colors} that meet all three criteria are: point A (x$_1$=21,812,415 MWh, x$_2$=77,102,051 MWh), point B (x$_1$=24,862,479 MWh, x$_2$=3,900,512 MWh), point C (x$_1$=45,223,754 MWh, x$_2$=74,516,399 MWh), and point D (x$_1$=47,475,469, x$_2$=0). Point B is the result we got from section 5. Also, notice that point C is definitely overproducing, so it is not considered in the calculations below. 
Next, we evaluate both the objective function from the previous sections and the alternative objective function at these corner points (A, B, and D):\\
257,326,658

\begin{itemize}
    \item For C(x$_1$, x$_2$) = 10.35x$_1$ + 19.51x$_2$, we get the results shown in Table \ref{tab:results_new_objective_function}.
    {
    \begin{table}[ht]
    \centering
    \begin{tabular}{|p{1.5cm}|p{4cm}|p{3cm}|}
    \hline
        Point & (x$_1$, x$_2$) & objective function value \\ \hline 
           A & (21812415, 77102051) &  1,730,019,510 $\$$\\ 
           B & (24862479, 3900512)  &  333,464,655 $\$$ \\ 
           D & (47475469, 0)  &  491,371,104 $\$$ \\ \hline
      
    \end{tabular}
    \caption{Values of the new objective function (\protect\ref{new_objective_function}) for points A, B, and D.}
    \label{tab:results_new_objective_function}
    \end{table}
    }
    
    \item For C(x$_1$, x$_2$) = 37.80x$_1$ + 58.62x$_2$, we get the results shown in Table \ref{tab:results_old_objective_function}.\\
    {
    \begin{table}[ht]
    \centering
    \begin{tabular}{|p{1.5cm}|p{4cm}|p{3cm}|}
    \hline
        Point & (x$_1$, x$_2$) & objective function value \\ \hline 
            A & (21812415, 77102051) &  5,344,231,517 $\$$\\ 
            B & (24862479, 3900512)  &  1,168,449,731 $\$$ \\
            D & (47475469, 0)  &  1,794,572,728 $\$$ \\ \hline
    \end{tabular}
    \caption{Values of the old objective function (\ref{old_objective_function}) for points A, B, and D.}
    \label{tab:results_old_objective_function}
    \end{table}
    }
    
\end{itemize}

Notice that the corner point that minimizes both objective function is point B. This is the result of wind energy being both cheaper to maintain and operate, and to build. So, even when we change the objective function, we still get the same suggested output.\\

\subsection{Recycling the materials}
The two main sources of renewable energy considered in this paper are wind and solar energy. Wind energy is obtained using wind turbines and solar energy with solar panels. Both of these devices have been originally produced from non-recycled materials, and that has become a problem for the environment. Since the lifespan of a solar panel is approximately 30-35 years the lifespan of a wind turbine is around 20 years, there is a growing need for technologies that will enable us to recycle both of them. Although most of the materials can be recycled, most of them have to be reused elsewhere and cannot be integrated into new wind blades or solar panels. 

\subsubsection{Wind blades}

It is estimated that 10kg of blade material is needed for every 1kW of new capacity. For our linear program, we consider 3.5 MW turbines, so each requires 35 tones of blade material. Since their lifespan is about 20 years, there is a need for recycling these blades. There are three main recycling technologies: mechanical recycling, thermal recycling, and chemical recycling. Although 99$\%$ of the materials can be reused, they cannot be used to make new blades \cite{WindRecycle}.

\subsubsection{Solar energy systems} \label{minimizing_emission_constraint}

There are three kinds of solar panels, but the crystalline silicon (monocrystalline or multi-crystalline) have higher conversion efficiency than others, so they are presently the most widely used commercial solar panels. After disassembly and extraction of the elements that compose such a panel, the weights of the various materials used are: glass 57.4$\%$, aluminum 12.7$\%$, adhesive sealant 10$\%$, silicon 3.1$\%$, and other 19.5$\%$. Some of these materials can be reused for new solar panels, but others lose their qualities and have to be reused elsewhere. Since solar energy is relatively new, most panels have not yet reached their end-of-life, and there is not much research on recycling technologies yet \cite{SolarRecycling}.

\section{Minimizing emissions instead of costs.} \label{minimizing_emissions_constraint}
In this study, the financial limitation of the project was considered of primary importance, and the goal was to minimize costs while keeping emissions low. However, an alternative model is presented in this section of the Appendix, where we minimize emissions while keeping the costs constrained. We chose to alter the model \ref{section_6_linear_model_program} from section \ref{section_6_space_solar}. Our objective function is now minimizing emissions (\ref{new_objective_emissions}), so there is no constraint on them. To keep costs low, we change the budget constraint coefficient to those used in function (\ref{old_objective_function}), which is represented in Inequality (\ref{budget_cost}. Here is the resulting model:

\quad Minimize 
\begin{equation}
    C(x_1, x_2) =  4970 x_1 + 45,000 x_2 
    \label{new_objective_emissions}
\end{equation} 
\quad Subject to:
\begin{equation}
    0.3760 \times x_1 + 0.01 \times x_2 \geq 7.069 \times 10^6 \text{MWh}
\end{equation}
\begin{equation}
    0.3775 \times x_1 + 0.9797 \times x_2 \geq 13.192 \times 10^6 \text{MWh}
\end{equation}
\begin{equation}
    0.2456 \times x_1 + 0.01 \times x_2 \geq 6.006 \times 10^6 \text{MWh}
\end{equation}
\begin{equation}
    37.80 x_1 + 58.62 x_2 \leq \$ 2 \times  10^9 
\label{budget_cost}
\end{equation}
\begin{equation}
    1065.6 x_1 + 204.5 (x_2 - 2,190,438 \text{MWh}) \leq 50,589,860,000 \text{ ft$^2$}
\end{equation}
\begin{equation}
    0 \leq x_1 , x_2
\end{equation}\\

\subsection{Results}
Solving this linear system, we get exclusively wind energy. This is expected, since wind is the cleanest energy source by far. It also manages to satisfy all constraints, so there is no need to add solar energy into the solution, as a less environment friendly option. 

Given this advantage that wind energy has over the rest of sources we consider, the rest of the models resulted in either exclusive use of wind energy or, in case of the model in section \ref{section_7_nuclear}, a combination of wind and nuclear because of the limitations on land-use.

\section{Energy systems of choice}
\subsection{Solar systems}

There are two main types of technologies for converting solar energy into electricity that are used at large scales: solar thermal technologies and photovoltaic systems. In this project, we will be considering the photovoltaic (PV) power plants, because they offer several advantages \cite{PV_vs_CSP}: 

\begin{itemize}
    \item It can use both the direct and diffused component of solar
radiation.
    \item It is suitable in areas with low direct radiation.
\end{itemize}

However, one of the major trade-offs of these solar systems is that their emissions are higher than for the thermal technologies. The emissions are still much lower than for conventional non-renewable sources, but are higher than some alternatives of solar systems.

\subsection{Wind turbines}
The size and height of turbines has been increasing every year. A larger turbine will produce more energy than a smaller one, so less turbines are needed to produce the same amount. A logical conclusion is: we should use the largest turbines we have. However, this is not as simple as might seem to be. Larger turbines are more difficult to transport, since they cannot be folded once they have been constructed. 

The average capacity of all the newly installed turbines in the US is 3 MW. Hence, this will be the theoretical power we are going to use for our calculations in this project. We will also be assuming that our turbines are at least 100 meters (330 feet) tall, since the speed of the wind increases at higher altitudes \cite{Wind_turbines_of_choice}. 

Illinois and the Midwest in general are considered to be great areas to build wind power plants since these are two regions with higher-than-average wind shear. Additionally, there are several manufacturing facilities located in Illinois and adjacent states, which allows us to assume that the transportation costs will be at most average, if not lower. Refer to the two maps below: \\

\includegraphics[width=6cm]{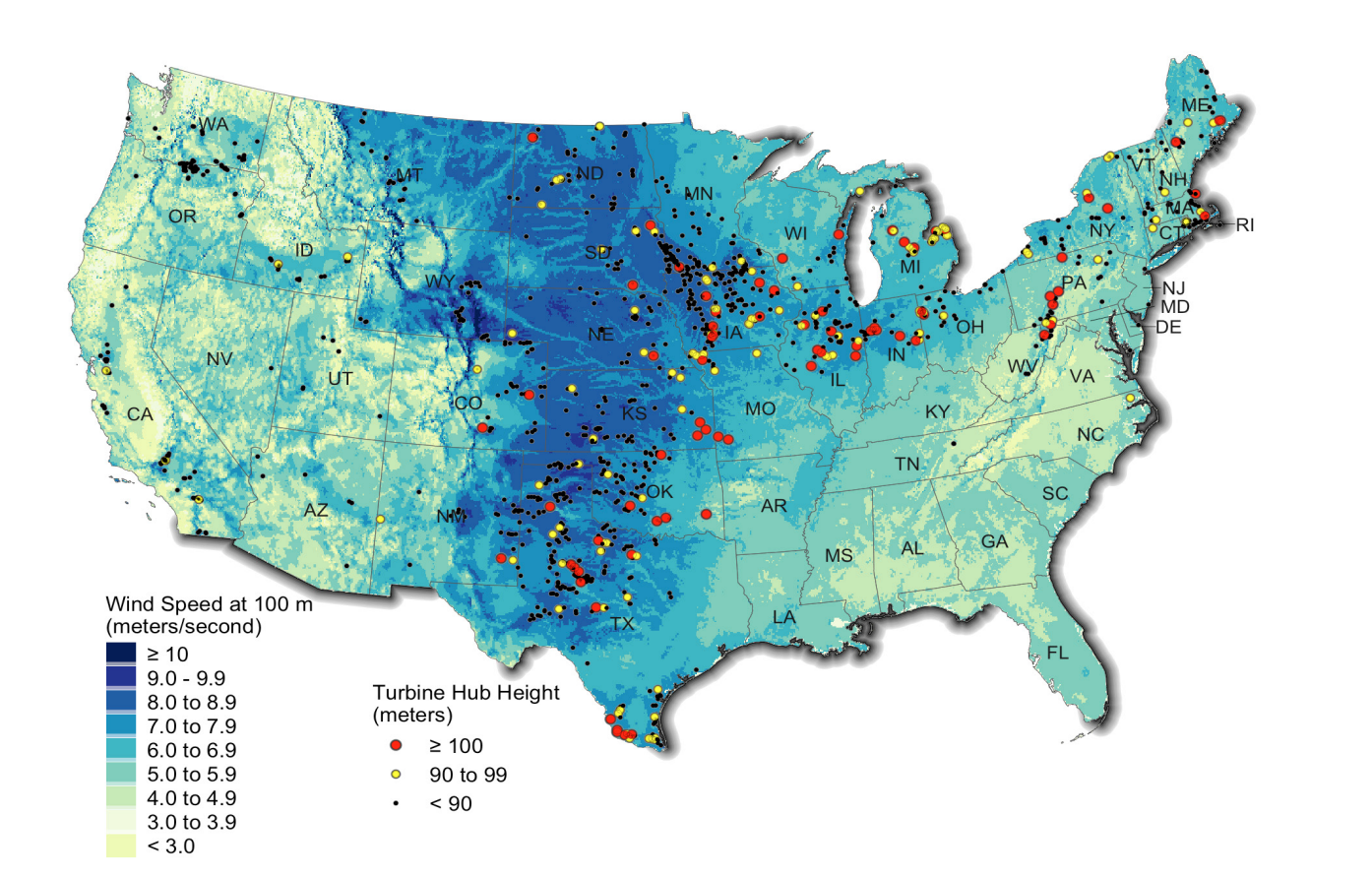}
\includegraphics[width=6cm]{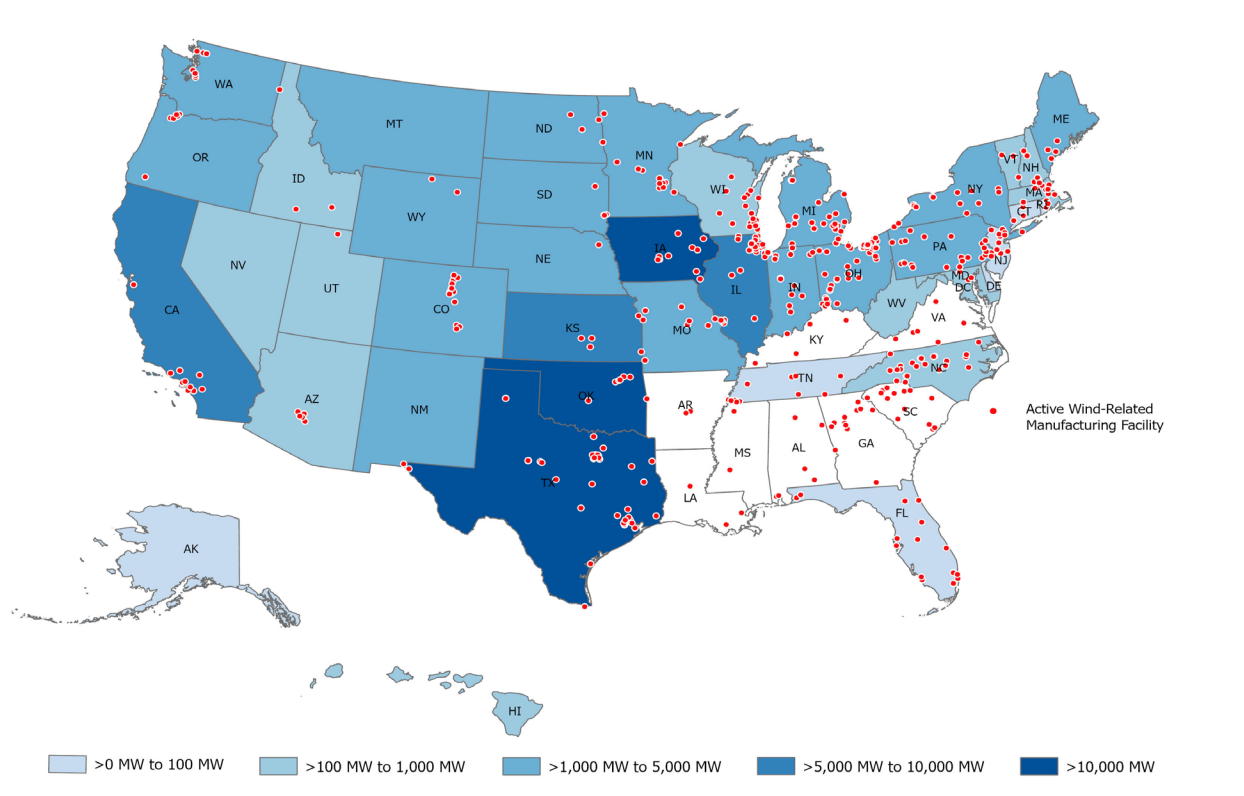}\\

From the first map, we can see that Illinois and the Midwest have rather high wind speed at 100m. From the second map, we can see the locations of several wind turbine manufacturing facilities \cite{Wind_manufacturing_speed}.

\clearpage

\end{document}